%Basic AmsTeX article template

\input amstex

\overfullrule=0pt
\define\ssk{\smallskip}
\define\msk{\medskip}

\define\del{$\partial$}
\define\ctln{\centerline}
\define\nidt{\noindent}
\define\blkbx{{\vrule height1.5mm width1.5mm depth1mm}}

\documentstyle{amsppt}

\topmatter

\title\nofrills Free Seifert surfaces and disk decompositions\endtitle
\author  Mark Brittenham \endauthor

\leftheadtext\nofrills{Mark Brittenham}
\rightheadtext\nofrills{Free Seifert surfaces and disk decompositions}

\affil   University of  North Texas \endaffil
\address   Department of Mathematics, University of North Texas, 
Denton, TX 76203  \endaddress
\curraddr Department of Mathematics, University of Nebraska - Lincoln,
Lincoln, NE 68588-0323 \endcurraddr
\email mbritten\@math.unl.edu \endemail
%\dedicatory       \enddedicatory
%\date                    \enddate 
\thanks   Research supported in part by NSF grant \# DMS$-$9704811 \endthanks
%\translator         \endtranslator
%\subjclass            \endsubjclass    

\abstract In this paper we construct families of 
knots which have genus one free Seifert surfaces which 
are not disk decomposable. \endabstract

\endtopmatter

\heading{\S 0 \\ Introduction}\endheading

A Seifert surface for a knot $K$ in the 3-sphere is an embedded
 orientable surface $\Sigma$, whose boundary equals the knot $K$.
Equivalently, it is a properly embedded orientable surface $\Sigma$ in 
the exterior $X(K)$ of $K$, whose boundary equals the longitude of $K$.
A Seifert surface $\Sigma$ is {\it free} if $\pi_1(S^3\setminus\Sigma$) is
a free group; equivalently, $S^3\setminus$int$N(\Sigma)$ is a handlebody.

Seifert's algorithm [Se] will always build a free Seifert surface 
for a knot $K$.
In [Br] we showed that not all free Seifert surfaces can be built by Seifert's
algorithm; we exhibited a family of hyperbolic
knots having {\it free genus} one, whose
surfaces built via Seifert's algorithm must always have large genus.

In so doing, we introduced a fairly general procedure for producing knots 
with genus
1 free Seifert surfaces. In this paper we show that many of these 
surfaces fail to
be disk decomposable.

A sutured manifold ($M$,$\gamma$) is a compact 3-manifold $M$ together
with a collection of disjoint embedded loops $\gamma$ in \del$M$, called
the sutures.
(Since we will apply this theory to knots and Seifert surfaces, we will 
suppress
the possibility that whole components of \del$M$ are sutures). The 
boundary of $M$ can
be expressed as \del$M$ = $R_+(\gamma)\cup R_-(\gamma)$, with 
$R_+(\gamma)\cap R_-(\gamma)$ = $\gamma$. We give $R_+$ a transverse
orientation pointing into $M$, and $R_-$ a transverse orientation
pointing out of $M$. We think of each component of
 $\gamma$ as having as having a transverse orientation pointing from its
$R_+$ side to its $R_-$ side, . For further details, see [Ga1]. 

A {\it decomposing surface} for ($M$,$\gamma$) is a properly embedded, 
transversely oriented surface $F$ which is transverse to 
$\gamma$. We can , by matching the transverse orientations for 
$R_+$, $R_-$, and $F$, endow $M$ split open along $F$, 
$M|F$, with the structure of a sutured manifold; see Figure 1.
The new sutures for $M|F$ are obtained as an `oriented sum'
of $\gamma$ and \del$F$. A sequence of such splittings is called
a {\it sutured manifold decomposition} of ($M$,$\gamma$).

\ssk

\input epsf.tex

\leavevmode

\epsfxsize=3.5in
\centerline{{\epsfbox{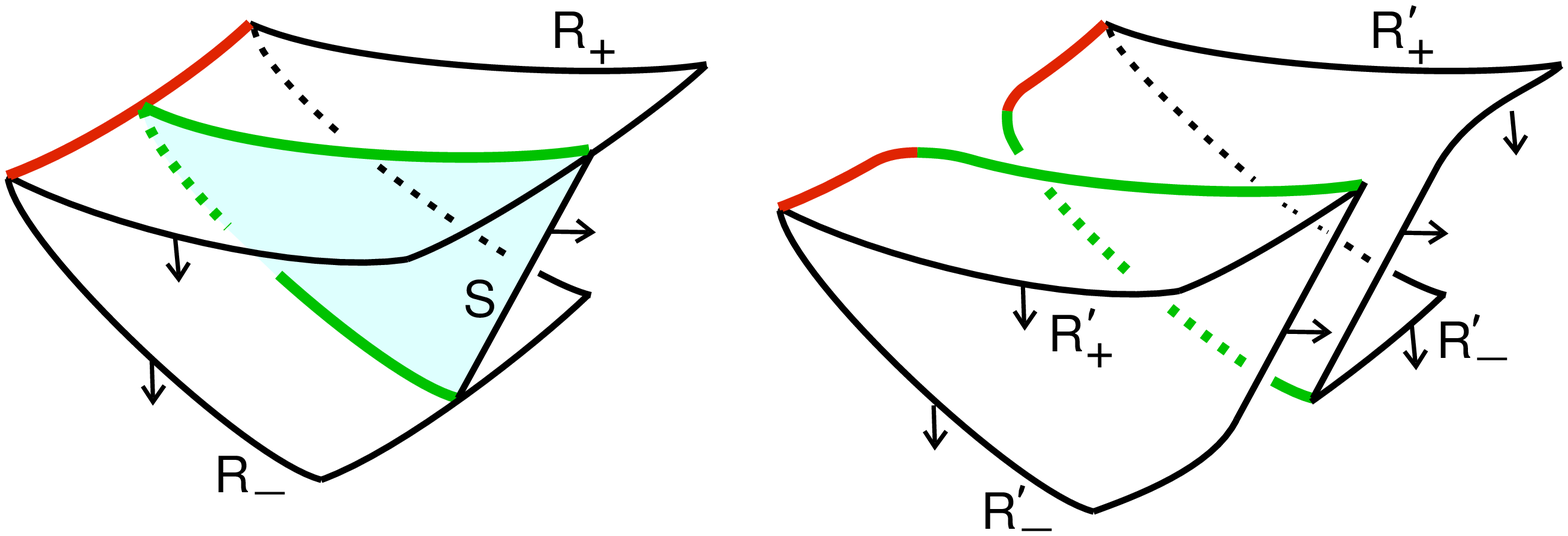}}}

\centerline{Figure 1}

\ssk

A Seifert surface $\Sigma$ is {\it disk decomposable} if the sutured manifold

\ctln{($S^3\setminus$int$N(\Sigma)$,$\Sigma\cap$\del$N(K)$) = 
($X_F,\gamma$)}

\nidt admits a sutured manifold decomposition whose decomposing 
surfaces are
all disks, ending with a sutured manifold which is the disjoint union
of sutured manifolds of the form ($B^3,e$), where $e$ is the equatorial 
circle of the 3-ball
$B^3$. (A posteriori, $S^3\setminus$int$N(\Sigma)$ is a handlebody, since
it
may be cut open along disks to 3-balls,
so $\Sigma$ is free.) By Gabai [Ga2], if $\Sigma$ is disk decomposable, then
the corresponding sutured manifold is {\it taut}, and so in particular 
$\Sigma$ has minimal 
genus among all surfaces representing its homology class. In other words, the
genus of $\Sigma$ equals the genus of $K$.

Disk decomposability therefore gives an effective way to compute the genus of a
knot. For example, Gabai [Ga2] has shown that every knot in the
standard tables [Ro] has a projection for which Seifert's algorithm
gives a disk decomposable surface.
A fairly natural natural question to ask, then, is: 
how can we tell, short of producing a set
of decomposing disks, that a Seifert surface is disk decomposable? Our main
result shows that being free and having minimal genus, which are necessary,
are not sufficient.

\proclaim{Theorem} There exist knots $K$ in $S^3$ which admit genus one 
incompressible free Seifert surfaces which are not disk decomposable.
\endproclaim

Our result leaves open the question of whether or not these knots admit {\it other} Seifert
surfaces which are disk decomposable; we discuss this possibility in the concluding section of 
the paper.

%%%%%%%%%%%%%%%%%%%%%%%%%%%%%%%%%%%%%%%%%%%%%%%%%%%%%%%%%%%%%%%%%%%%%%%%%%%%%%%%%%%%%%%%%%%%%%%%%%%%%

\ssk

\leavevmode

\epsfxsize=4.9in
\centerline{{\epsfbox{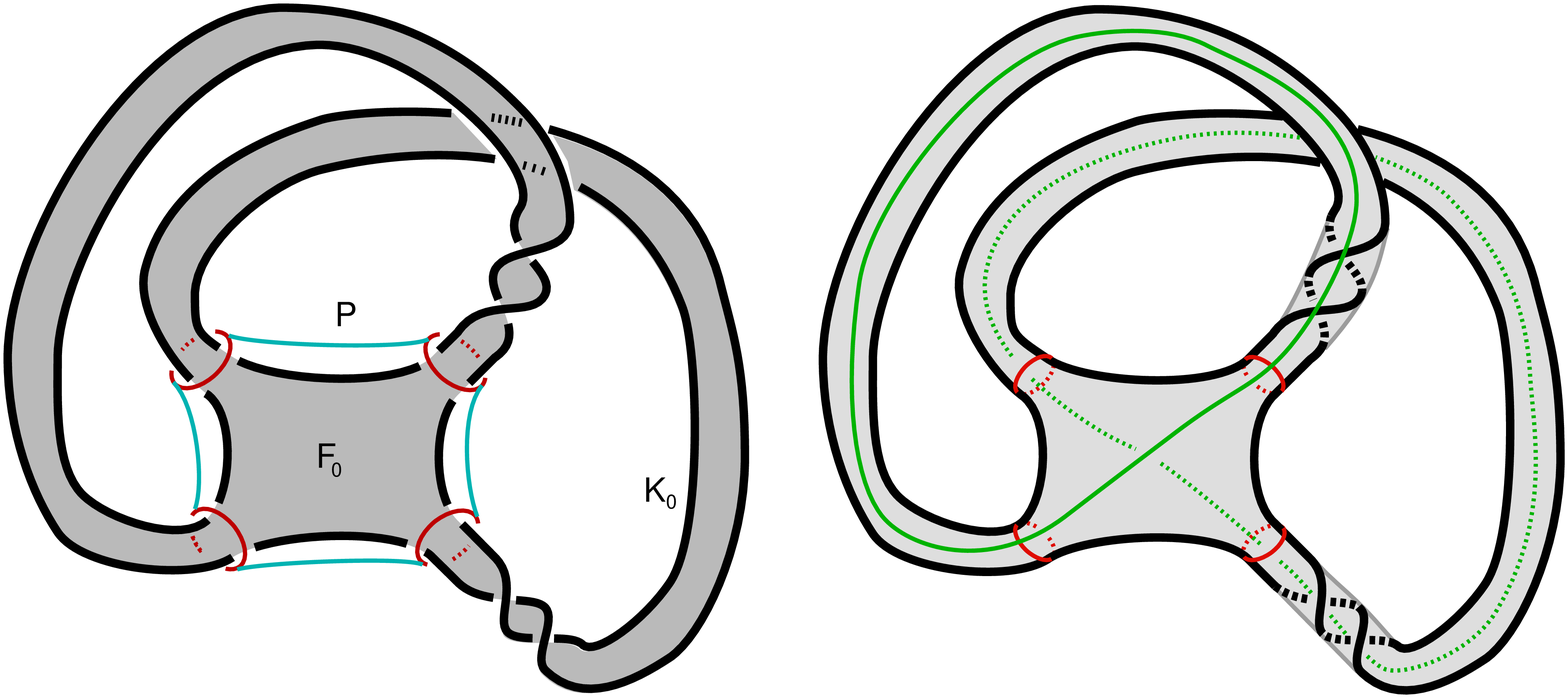}}}

\centerline{Figure 2}

\ssk

\heading{\S 1 \\ Building free Seifert surfaces}\endheading

In [Br] we showed that, for the knot $K_0$ and the free Seifert surface $F_0$ 
for $K_0$ in $S^3$, shown in Figure 2,
$1/n$ Dehn surgery on any loop $L$ in the 4-punctured sphere $P$ 
pictured there will essentially re-embed $K_0$ and $F_0$ as a new knot 
$K$ and free Seifert surface $F$ in $S^3$. 
(There is, in fact, nothing special about this knot; any free
Seifert surface for a knot admits similar 4-punctured spheres.)
What we will show now is that for appropriate
choices of $L$ and $n$, $F$ will be incompressible but not disk
decomposable. Our sutured manifold $X_F$ = $S^3\setminus$int$N(F)$ 
will be a genus-2
handlebody, and the suture will be a loop $\gamma$ = $F\cap$\del$X(K)$
which splits \del$X_F$ into two
once-punctured tori. The essential idea is that if $\gamma$ is 
complicated enough with respect to a set of cutting disks for $X_F$,
then $F$ must be incompressible in $X(K)$, but not disk decomposable. 

\ssk

\leavevmode

\epsfxsize=4.9in
\centerline{{\epsfbox{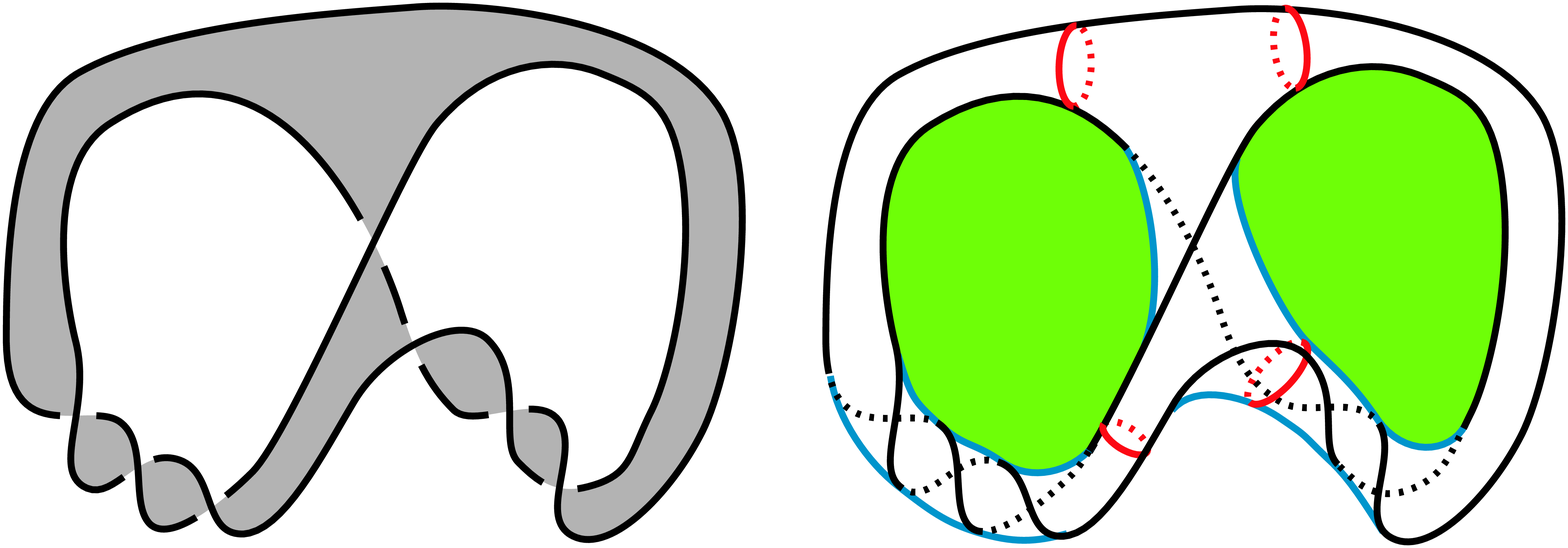}}}

\centerline{Figure 3}

\ssk

\leavevmode

\epsfxsize=4in
\centerline{{\epsfbox{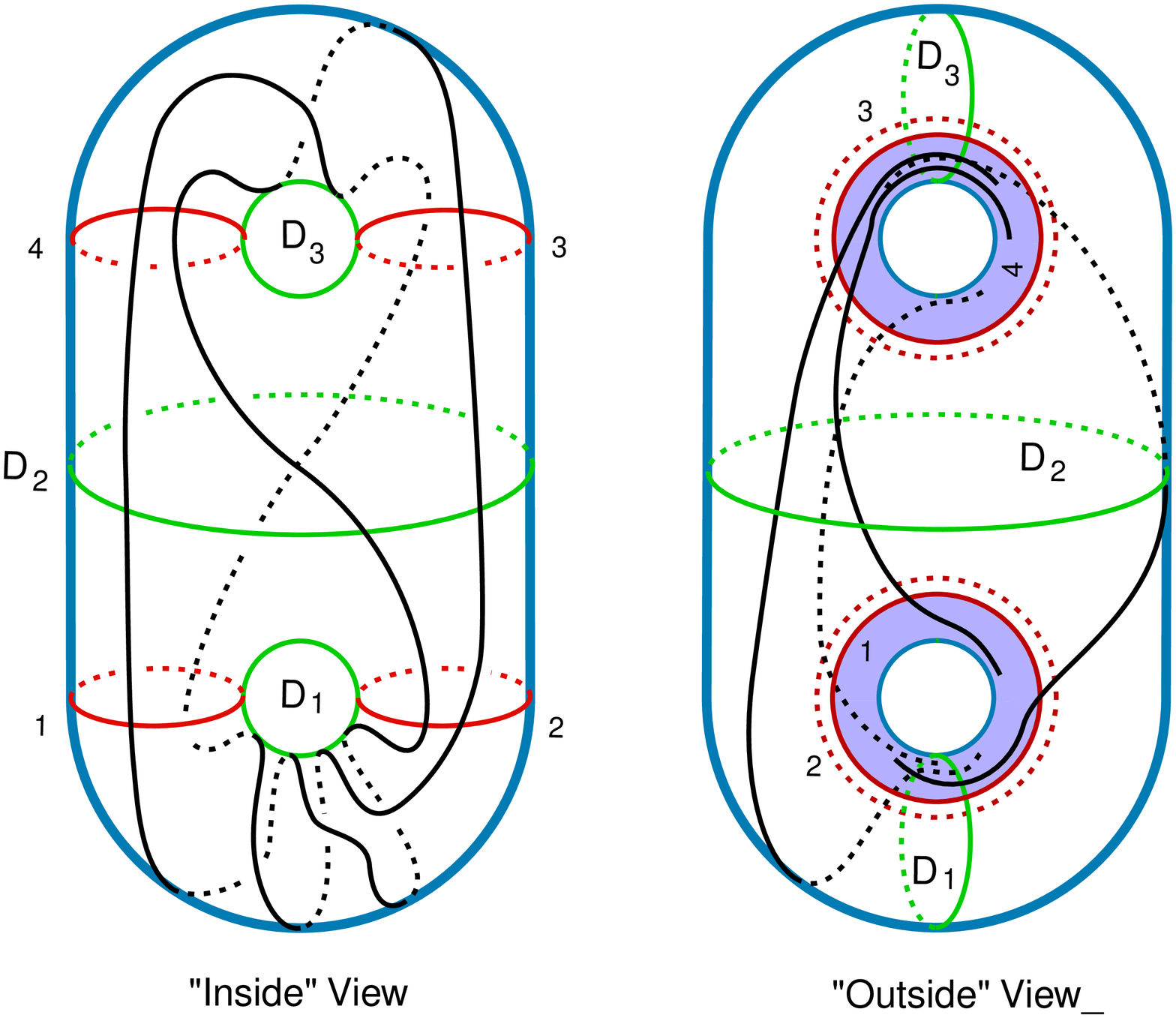}}}

\centerline{Figure 4}

\ssk

Our argument will be based on techniques of Goda [Go], who first showed that
there exist taut sutured handlebodies ($H$,$\gamma$) which are not 
disk decomposable. Our main task will be to show that his arguments can be 
applied to some of the sutured handlebodies built as in Figure 2. Because 
Goda's techniques use the standard view of a handlebody, as
the inside of a standardly embedded genus two surface, we first need to 
produce an `external' view of $X_{F_0}$. In other words, we need
to understand what our suture $\gamma$ looks like when $X_{F_0}$ is
pictured as the interior of a standardly embedded handlebody 
in $S^3$. This involves, 
essentially, determining what the two annuli 
in \del$X_{F_0}$ that are cut off by \del$P$ would look like on a standard 
handlebody, while keeping track of the pattern of intersections of
$\gamma$ and \del$P$ with a set of three `obvious' cutting disks for $X_{F_0}$, 
whose boundaries are shown in Figures 3 and 4. This pattern determines $\gamma$
up to homeomorphism, in fact, up to Dehn twists along the three cutting
disks, since these disks form a complete system of cutting disks for $X_{F_0}$.

This change of viewpoint is carried out in Figure 4. Our suture $\gamma$
can be thought of as four arcs lying on a 4-punctured sphere (essentially, $P$) in
$X_{F_0}$, together with two pairs of arcs spiralling through the complementary
annuli in \del$X_{F_0}$ . The amount of spiralling is determined by how many 
full twists we put in each arm of our original Seifert surface $F_0$, and will
not play a large role in our further discussions (although the 
{\it direction} of spiralling
{\it is} important).

\heading{\S 2 Choosing loops $L$\\ }\endheading

It is easy to see what effect $1/n$ Dehn filling on a loop $L$ in $P$ will have
on the picture of our sutured handlebody ($X_{F_0},\gamma$) above. The loop $L$ 
will lie
slightly inside of the 4-punctured sphere $P$ lying on \del$X_{F_0}$, and the 
disk $D$ it bounds will lie for the most part outside of the handlebody $X_F$, since 
it mostly lies in $N(F_0)$. Inside of our handlebody we will see only an annulus
running from $L$ to a parallel loop ($L^\prime$, say) in $P$. 
$1/n$ Dehn surgery along 
$L$ will have the effect of replacing the suture $\gamma$ with its
`sum` with $n$ parallel oriented copies of the loop $L^\prime$ (Figure 5). 
It will in fact be the result of applying $n$ Dehn twists in an annular neighborhood
of $L^\prime$ to $\gamma$.
This gives us a large family of sutured handlebodies to work with, each of which is 
realized as the complement of some genus one free Seifert surface in $S^3$. 

\ssk

\leavevmode

\epsfxsize=5in
\centerline{{\epsfbox{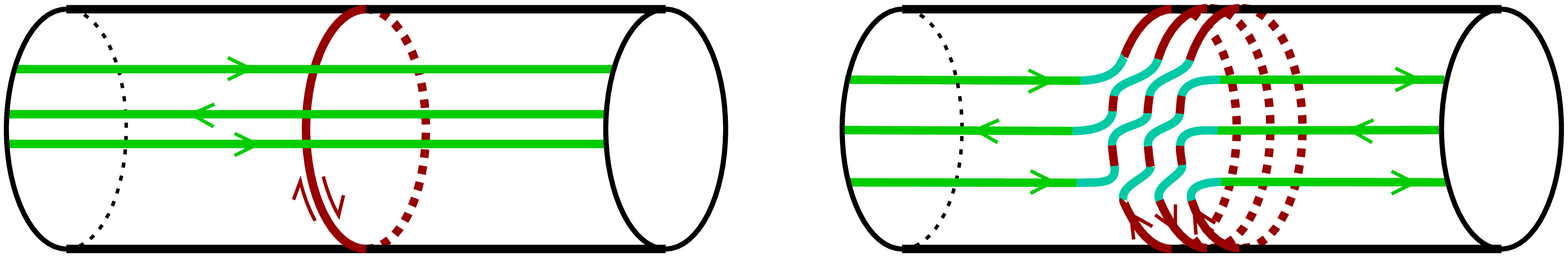}}}

\centerline{Figure 5}

We illustrate this with a somewhat more complicated loop $L$, in Figure 6; 
it meets $K_0$ in 
8 points, and so $n$ Dehn twists along $L$ will result in a knot $K$,
which on \del$H$ will be represented by the `sum' of $K_0$ and 8n parallel copies
of $L$. We show the results of one Dehn twist, in Figure 7. 

This new suture $C$ meets the standard cutting disks for the handlebody $H$ 
(which are
the three disks where a horizontal 
plane perpendicular to the paper meets the middle of the figure)
only in arcs joining distinct 
cutting disks. These arcs run, in each of the pairs of pants in
\del$H$, above and below the cutting disks, between any pair of the cutting 
disks. It is also easy to see that there are no trivial arcs, running from a 
cutting disk to itself. This implies that \del$H\setminus C$ is incompressible
in $H$ ([St],[Ko1]), and $(H,C)$ is therefore taut. 

\ssk

\leavevmode

\epsfxsize=3.6in
\centerline{{\epsfbox{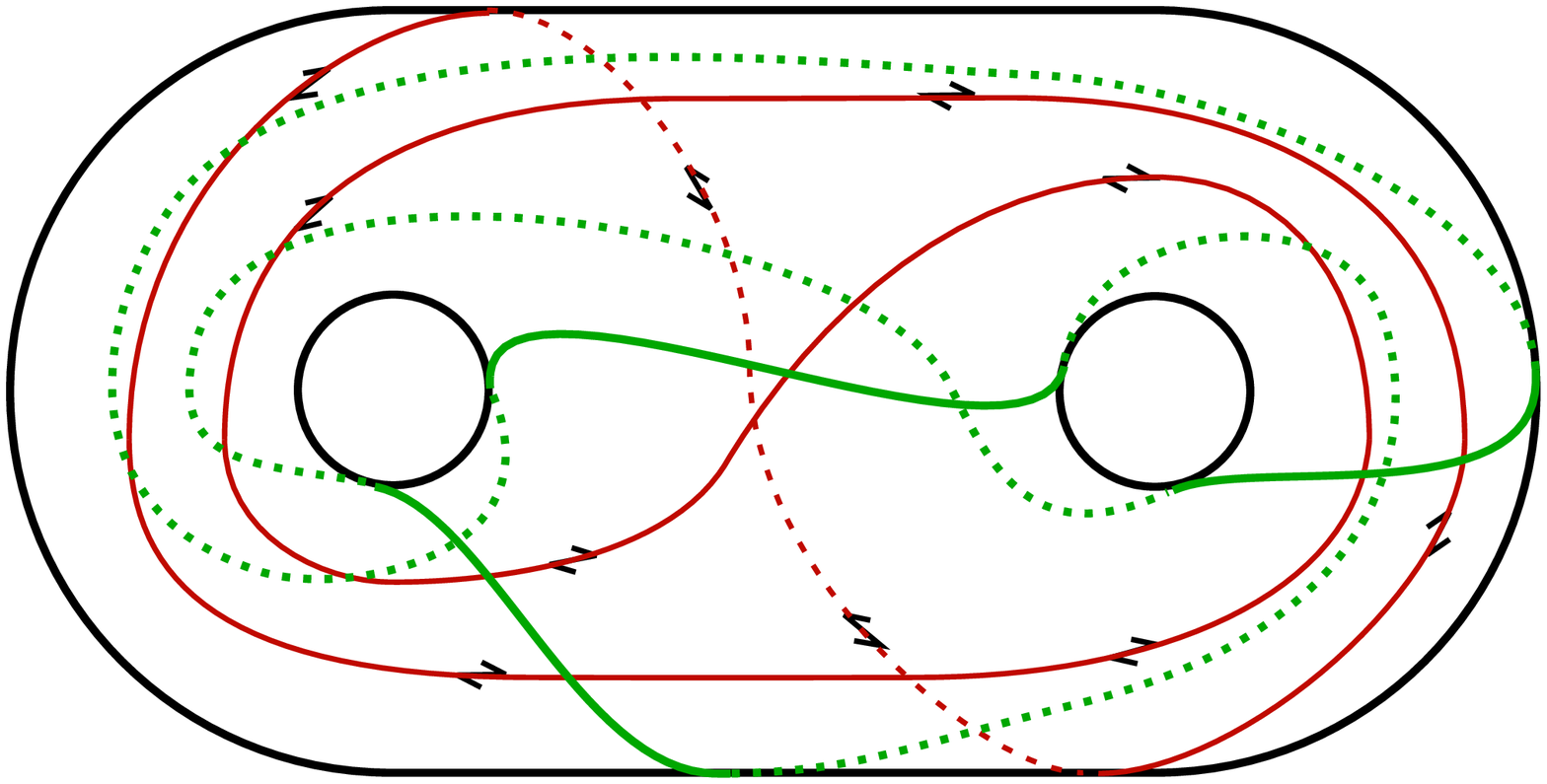}}}

\centerline{Figure 6}

\msk

\leavevmode

\epsfxsize=3.6in
\centerline{{\epsfbox{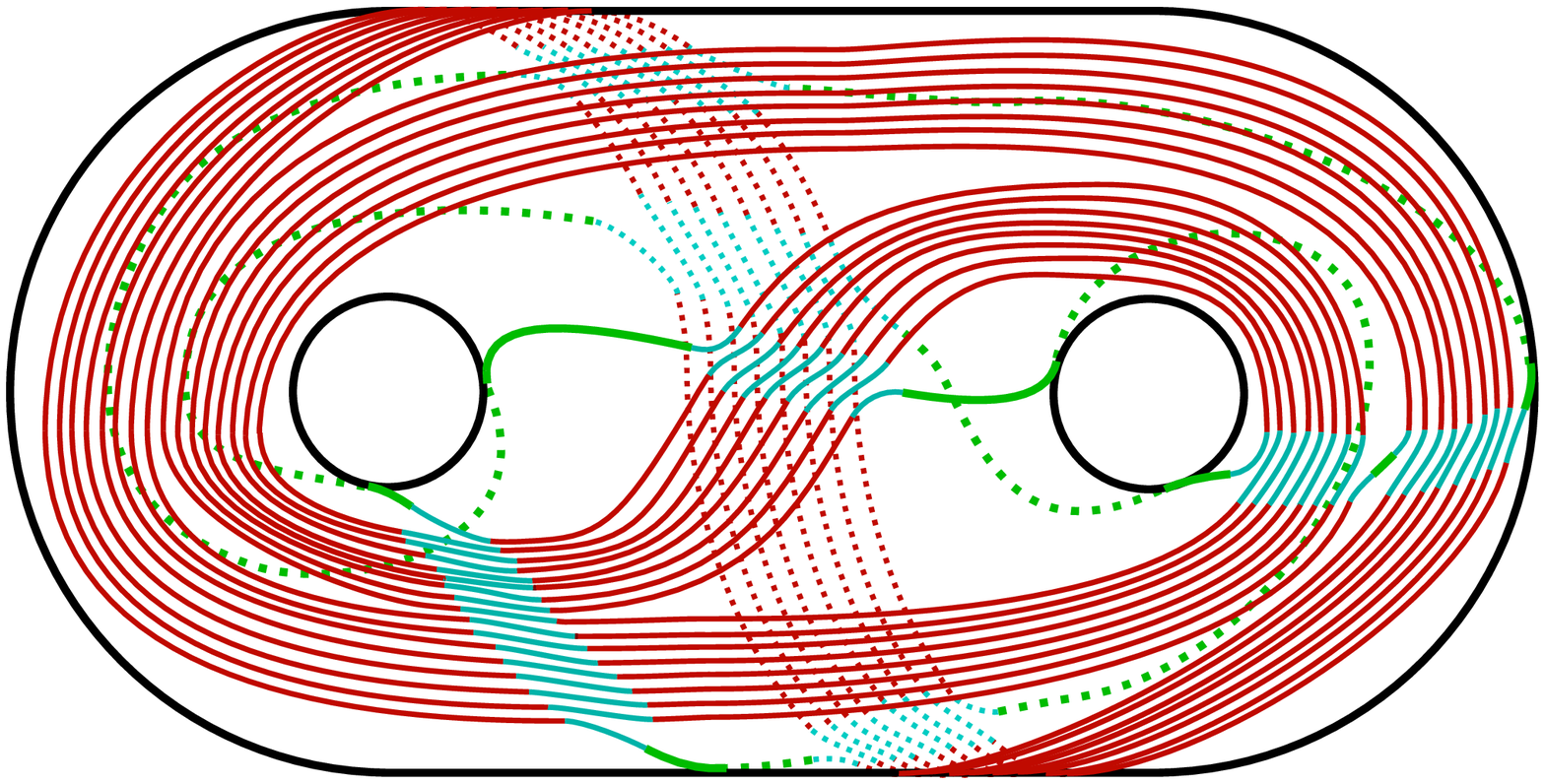}}}

\centerline{Figure 7}

\ssk

We can `encode' this construction, and the 
Dehn twisting information, into a train track $\tau$ on \del$H$ (Figure 8a) carrying
both $\gamma$, $L$, and the result $C$ of `right-handed' Dehn twists of $\gamma$ 
along $L$. This allows us 
to see, even for a large number of Dehn twists, that all of the loops so 
built represent sutures of taut sutured handlebodies, since it is easy
to see that any loop (which separates \del$H$) carried with full support by $\tau$ 
has arcs
running between any pair of the cutting disks $D_i$, on each side, as before, 
and has no trivial arcs. 
This is most easily seen by cutting \del$H$ (and $\tau$) open along our 
cutting disks (Figure 8b); the resulting train tracks carry no trivial arcs
or circles.

This curve $L$ (and the resulting sutures $C$) will, 
in the
end, still not be sufficiently `complicated' for our purposes. But several, 
which will be,
share many of the same properties, being carried by the same train track $\tau$.

\ssk

\leavevmode

\epsfxsize=3.5in
\centerline{{\epsfbox{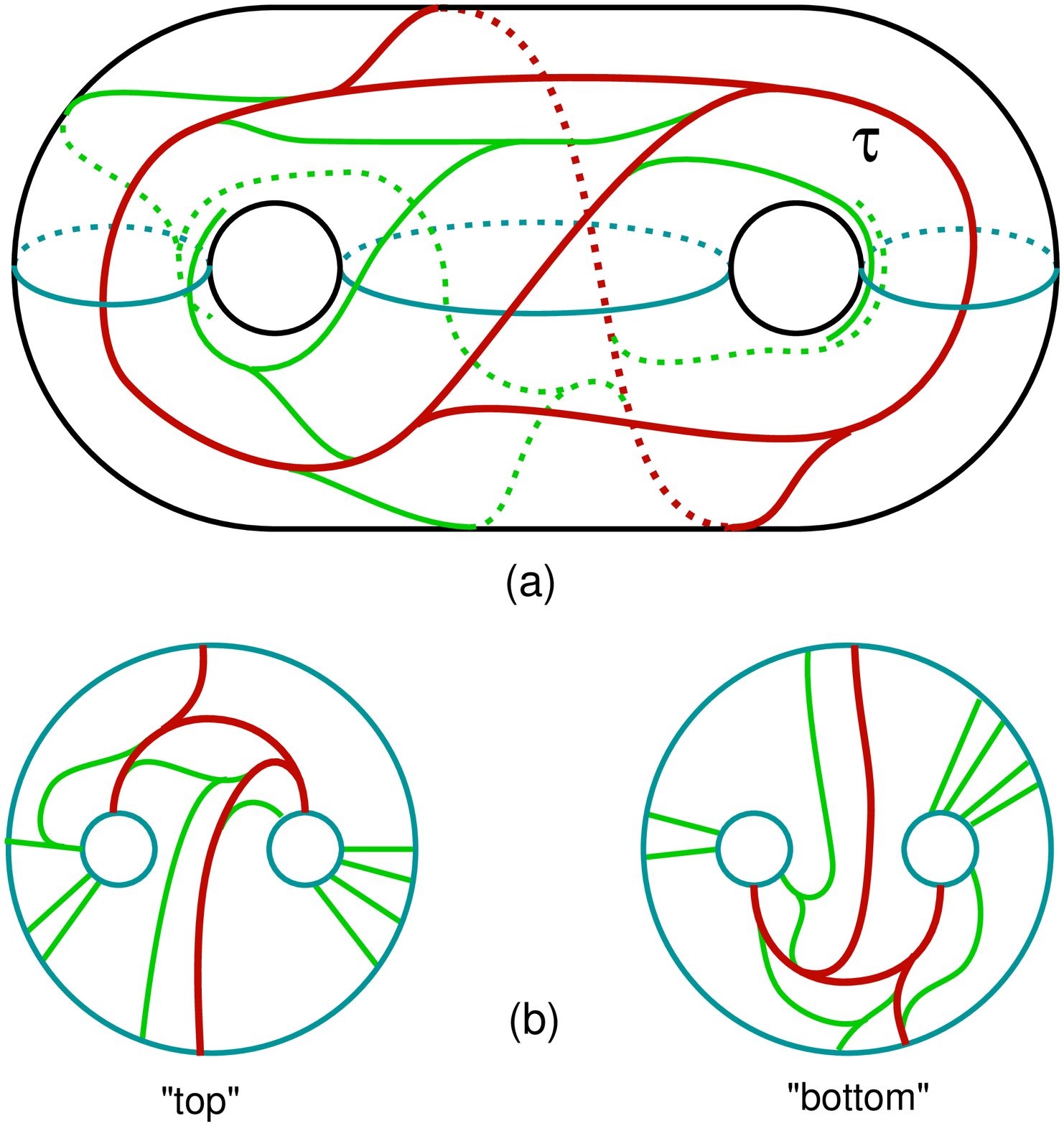}}}

\centerline{Figure 8}

\msk

\heading{\S 3 Complicated intersections imply no disk decomposition}\endheading

Goda [Go] determined sufficient conditions, based on the intersections of 
the suture $C$ in
the boundary of a genus two handlebody $H$ with a system of
cutting disks $D_1,D_2,D_3$ for $H$, to guarantee that the sutured handlebody 
$(H,C)$ is taut but not disk decomposable. 
We will prove here a 
slightly weaker form of Goda's criterion, which is sufficient for our
purposes. Note that any loop $C$ in \del$H$ locally separates \del$H$ (i.e., it
separates a neighborhood of itself), so we can always unambiguously talk 
about being on the `same side' of $C$ in \del$H$. 

\proclaim{Proposition 1} If $D$ is a compressing disk for \del$H$,
with \del$D$ transverse to $C$, such that 
$C|D$ contains three parallel arcs whose ends all lie on the same side of \del$D$,
then ($H,C$) is not disk decomposable along $D$.
\endproclaim

\ssk

\leavevmode

\epsfxsize=2in
\centerline{{\epsfbox{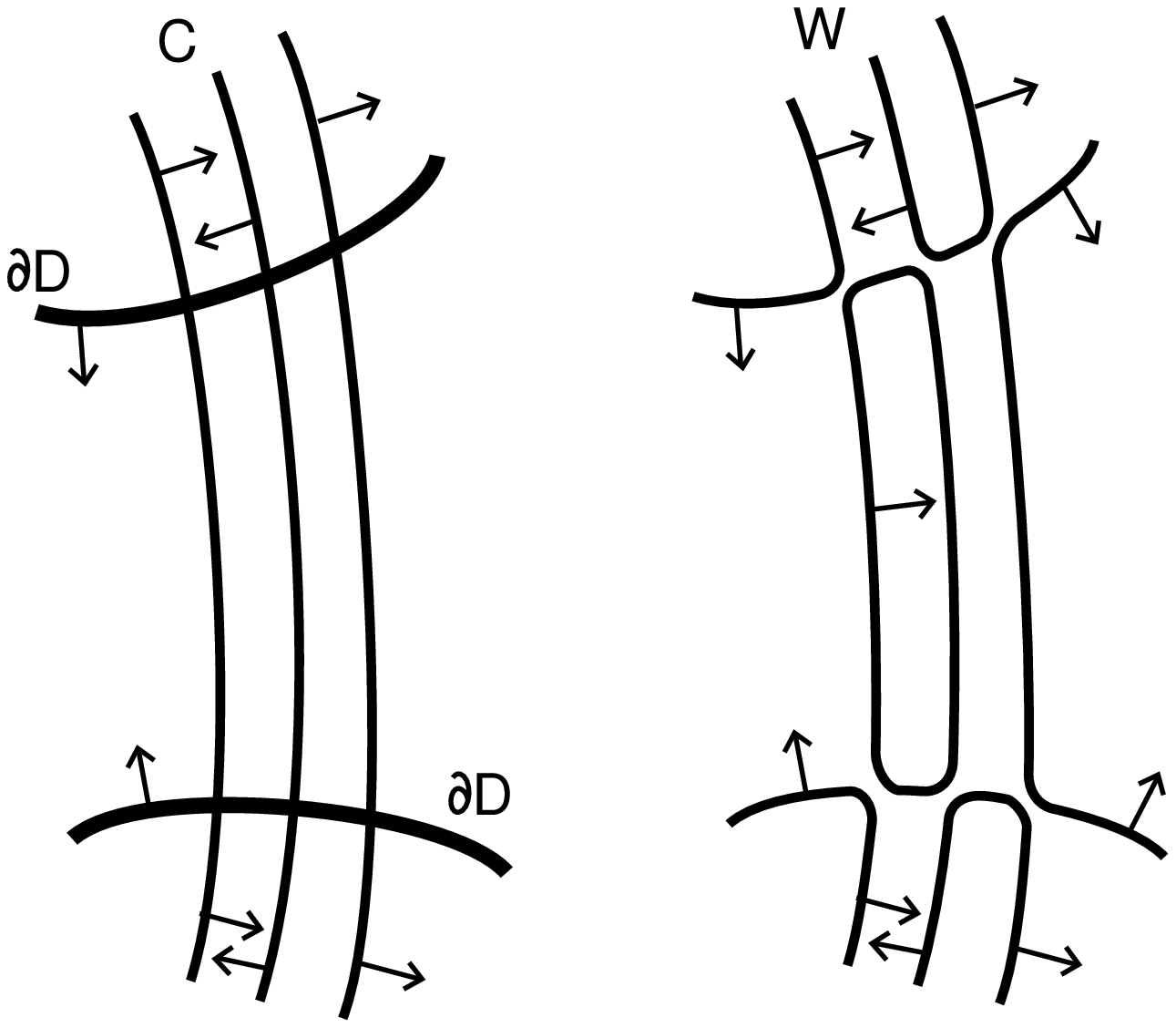}}}

\centerline{Figure 9}

\ssk

{\bf Proof:} The picture we have is as in Figure 9. (All other possible choices of 
normal orientation can be obtained from the one pictured by some combination of changing
every orientation or reflecting in a vertical axis, which will not change the essential
features of our argument.)
Given a transverse orientation on the disk $D$, 
the sutured manifold obtained by cutting $H$ along $D$ is one or two solid tori 
$H|D$, whose sutures are obtained by cutting and pasting
$C$ and \del$D$ near their points of intersection, as in Figure 1. 
However, because of our hypothesis, the sutures of the resulting sutured solid 
torus 
$(M,C^\prime)$ will include a component which is null-homotopic in \del$M$ 
(Figure 9), 
and so $(M,C^\prime)$ cannot be taut. The key point here is that 
$C$ separates \del$H$, and
so the transverse orientations of $C$, seen along \del$D$, must alternate.\blkbx

\msk

The reader can note that in Figure 7, each of the disks $D_i$ will have 3 such arcs on
each side. We illustrate one such collection in Figure 10. 

\ssk

\leavevmode

\epsfxsize=3.6in
\centerline{{\epsfbox{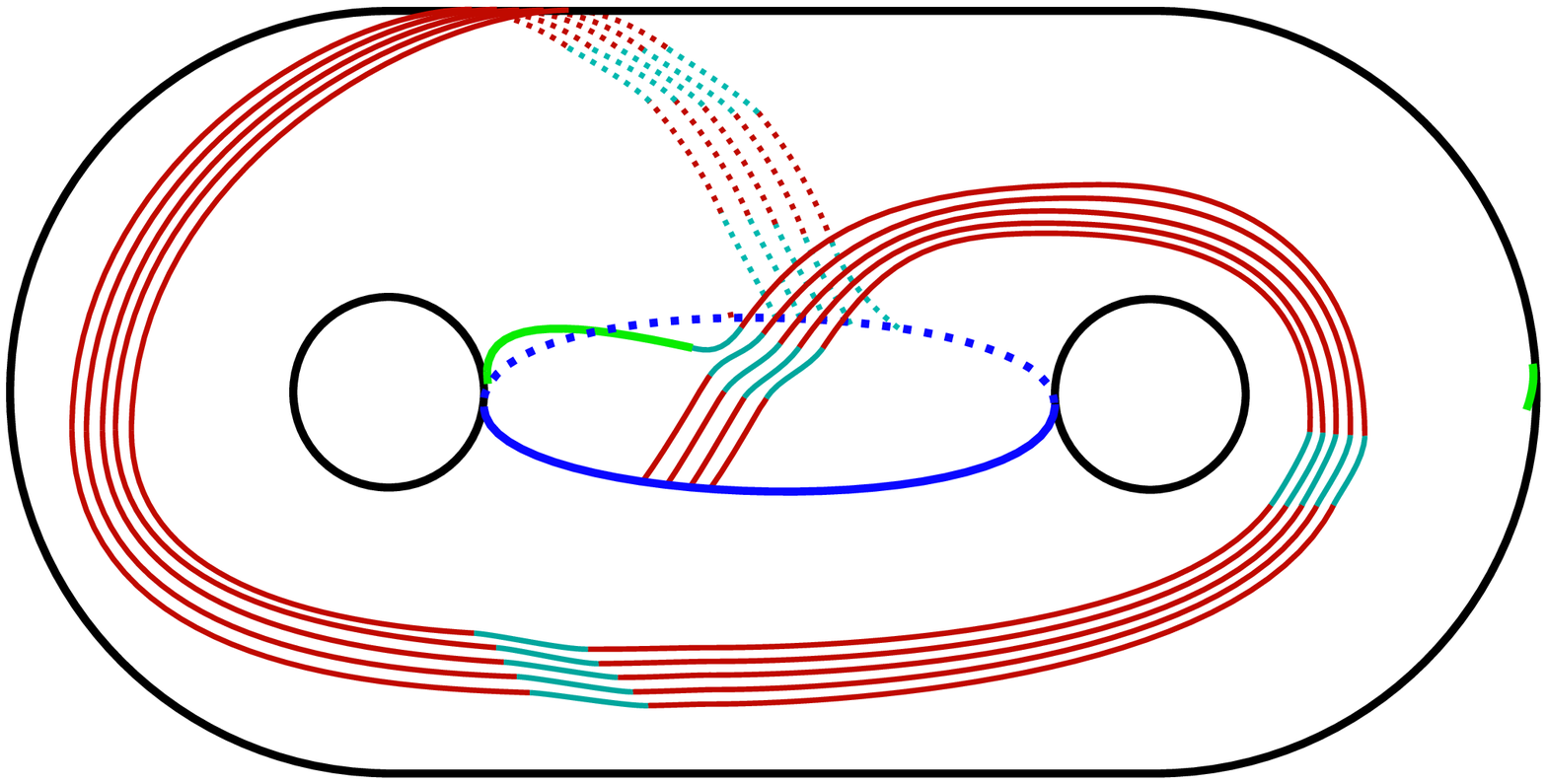}}}

\centerline{Figure 10}

\msk

We now assume that $C$ satisfies the conditions thus far introduced: the cutting disks
$D_i$ cut $C$ into arcs in the two pairs of pants 
\del$H|$(\del$D_1\cup$\del$D_2\cup$\del$D_3$)=$P_1\cup P_2$. Each arc joins distinct
\del-components in the $P_i$, and there are
arcs running between all possible pairs of \del-components of the $P_i$. We also have, 
for each disk
$D_i$, a set of three parallel arcs in \del$H|$\del$D_i$, as in Proposition 1.

\msk

\proclaim{Proposition 2} Any disk $D$, isotopic to one of the disks $D_1$, $i=1,2,3$ and
transverse to $C$, is not a decomposing disk for $(H,C)$.\endproclaim

{\bf Proof:} This is essentially Claim 3.6 of [Go]; for completeness, we reproduce the 
argument here, since many of the same ideas will be used later.

Without loss of generality, we may assume the $D$ is isotopic to $D_1$; then by 
[Ep,Lemma 2.5] there is an innermost disk $\Delta$ in \del$H$ whose boundary 
consists of an 
arc $\alpha$ of \del$D_1$ and an arc $\beta$ of \del$D$. $C$ intersects $\Delta$
in arcs, and by our hypothesis, none of these arcs have both endpoints on $\alpha$.
If any have both endpoints on $\beta$, then there is an outermost such arc $\delta$; 
but then it is easy to see that either decomposing ($H,C$) along $D$
yields a trivial suture, implying the $(H,C)$ is not disk decomposable along
$D$ (Figure 11a), or we may isotope $C$ across the outermost disk cut off by
$\delta$, without altering the sutured manifold obtained  by decomposing along $D$
(Figure 11b). Continuing, we can remove all such trivial intersections of $C$ with $D$
(or obtain our desired conclusion).

\ssk

\leavevmode

\epsfxsize=4in
\centerline{{\epsfbox{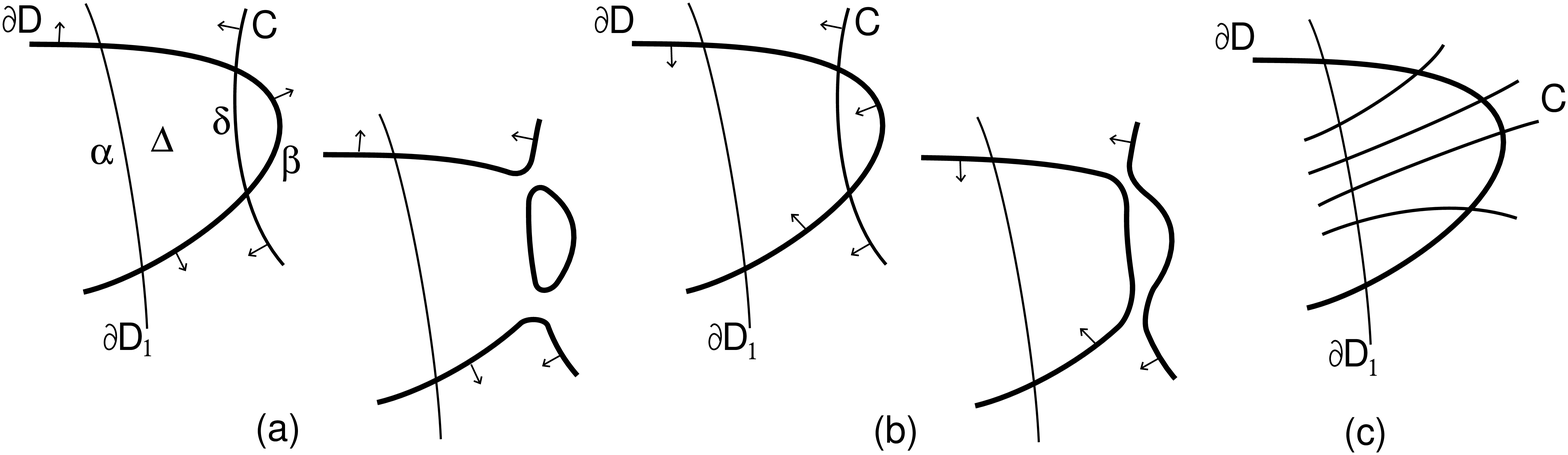}}}

\centerline{Figure 11}

\ssk

We may therefore assume that $C$ meets $\Delta$ only in arcs running from $\alpha$ 
to $\beta$, which must therefore all be parallel to one another (Figure 11c),
and so we can isotope \del$D$ across $\Delta$, removing two points of
intersection of \del$D$ with \del$D_1$, without changing the intersections of
\del$D$ with $C$. Continuing, we can then assume that \del$D$ and \del$D_1$
are disjoint, and so by [Ep,Lemma 2.4] they cobound an annulus $B$. By the same 
argument, we may assume that $C$ meets $B$ only in arcs running from \del$D$ to
\del$D_1$, and so we may isotope \del$D$ to \del$D_1$ without changing the
intersections of \del$D$ with $C$. Therefore the sutured manifold resulting from
decomposing along $D$ is identical with the one obtained by decomposing along 
$D_1$. But by our hypotheses and Proposition 1, $(H,C)$ is not
disk decomposable along $D_1$, and so it cannot be disk decomposable along $D$.
\blkbx

\ssk

Next we give a criterion which is sufficient to guarantee that every compressing
disk for \del$H$ has a trio of parallel arcs in $C$.

%%%%%%%%%%%%%%%%%%%%%%%%%%%%%

\proclaim{Proposition 3} 
Suppose that for every pair of cutting disks $D_i$ and $D_j$, $i\neq j$, and for each 
side of \del$D_i\subseteq$\del$H$, there is a collection of three parallel subarcs of 
$C$, with endpoints on the same side of \del$D_i$, which on both ends 
cross $D_j$ immediately before meeting $D_i$ (see Figure 12). Then for every 
disk $D$ in $H$, with \del$D\subseteq$\del$H$ transverse to $C$, $(H,C)$ is not
disk decomposable along $D$. \endproclaim

{\bf Proof:} Suppose that $D$ is a decomposing disk for $(H,C)$. By Proposition 2,
$D$ is not isotopic to any of the $D_i$. Because it is a compressing disk, \del$D$ cannot
be trivial in \del$H$. But since every simple loop in a pair of pants is either 
trivial or isotopic to one of the \del-components, this means that \del$D$ cannot be
isotoped to be disjoint from all of the \del$D_i$; it would then lie in one of our two
complementary pairs of pants.

\ssk

\leavevmode

\epsfxsize=2.4in
\centerline{{\epsfbox{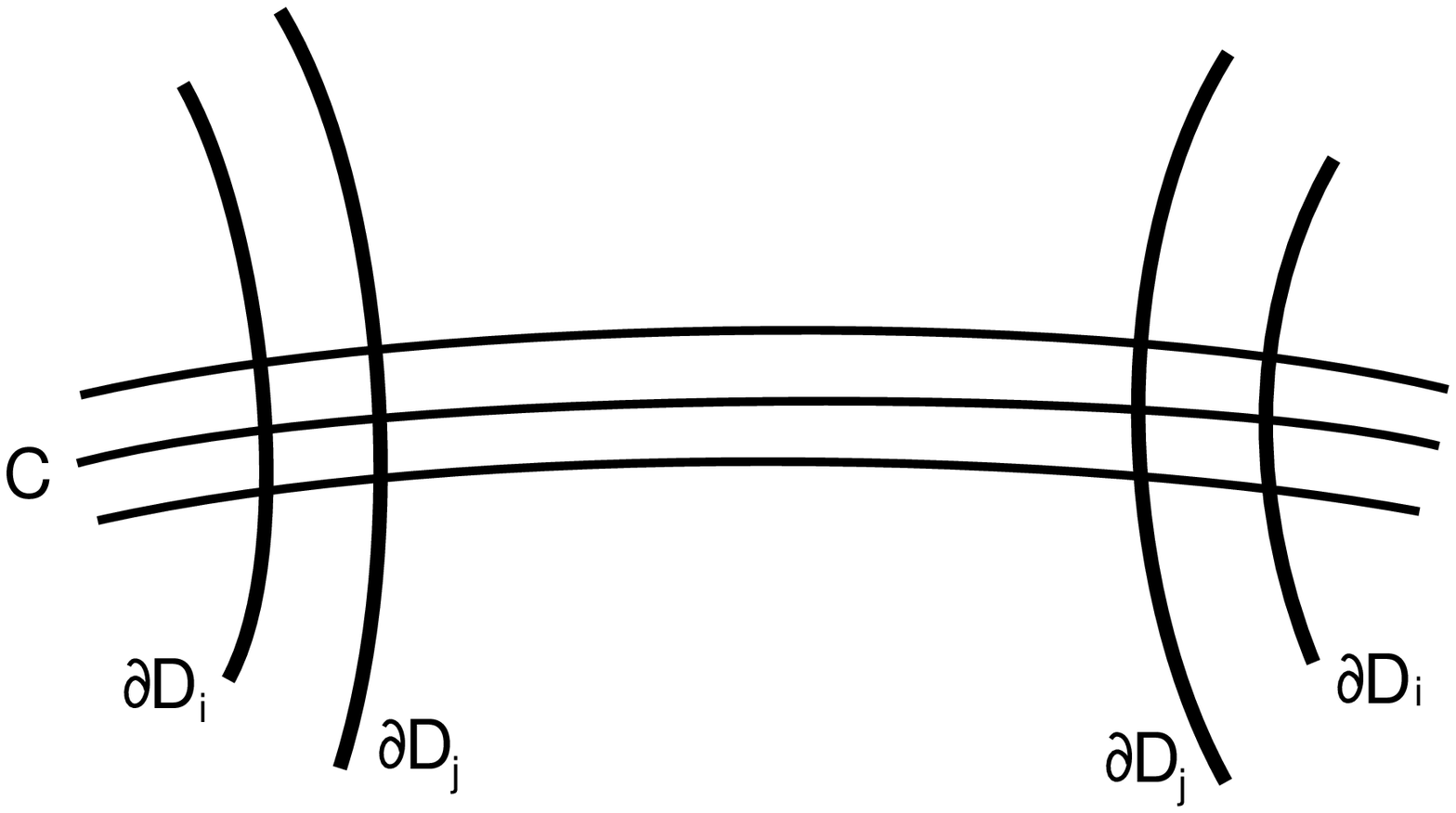}}}

\centerline{Figure 12}

\ssk

Consider an arc $\delta$ of $D\cap(D_1\cup D_2\cup D_3)\subseteq D$ which is outermost 
in $D$. The arc $\beta$ of \del$D$ which $\delta$ cuts off then lies in one of our two 
pairs of pants, call it $P$. If $\beta$ is a trivial arc in $P$, then, together 
with an arc $\alpha$
in one of the \del$D_k$, it bounds a disk $\Delta$ in
$P$. The suture $C$ meets $\Delta$ in arcs, and, by applying the 
argument of the previous proposition, we may assume that each arc runs
from $\alpha$ to $\beta$, since, if not, then either decomposing along 
$D$ will create a trivial
suture, implying that $D$ is not a decomposing disk for $H$, or we can 
isotope \del$D$ across $C$ without changing what the sutures in the
sutured manifold obtained by splitting along $D$ will look like.
The pictures are identical to those of Figures 11a and b. 
A trivial arc cannot lie on the $\delta$-side of $\Delta$, by hypothesis.

But then, as before, we can isotope 
\del$D$ across $\Delta$ to reduce the number of points
of intersection of \del$D$ with the \del$D_i$, without changing the sutured manifold
$(H|D,W)$. After repeatedly carrying out these isotopies, we can then assume
that every outermost arc in $D$ is non-trivial. By our argument above, there 
must be at least
one non-trivial arc, $\alpha$, since otherwise $D$, hence \del$D$, is
disjoint from the $D_i$. 

\ssk

\leavevmode

\epsfxsize=1.5in
\centerline{{\epsfbox{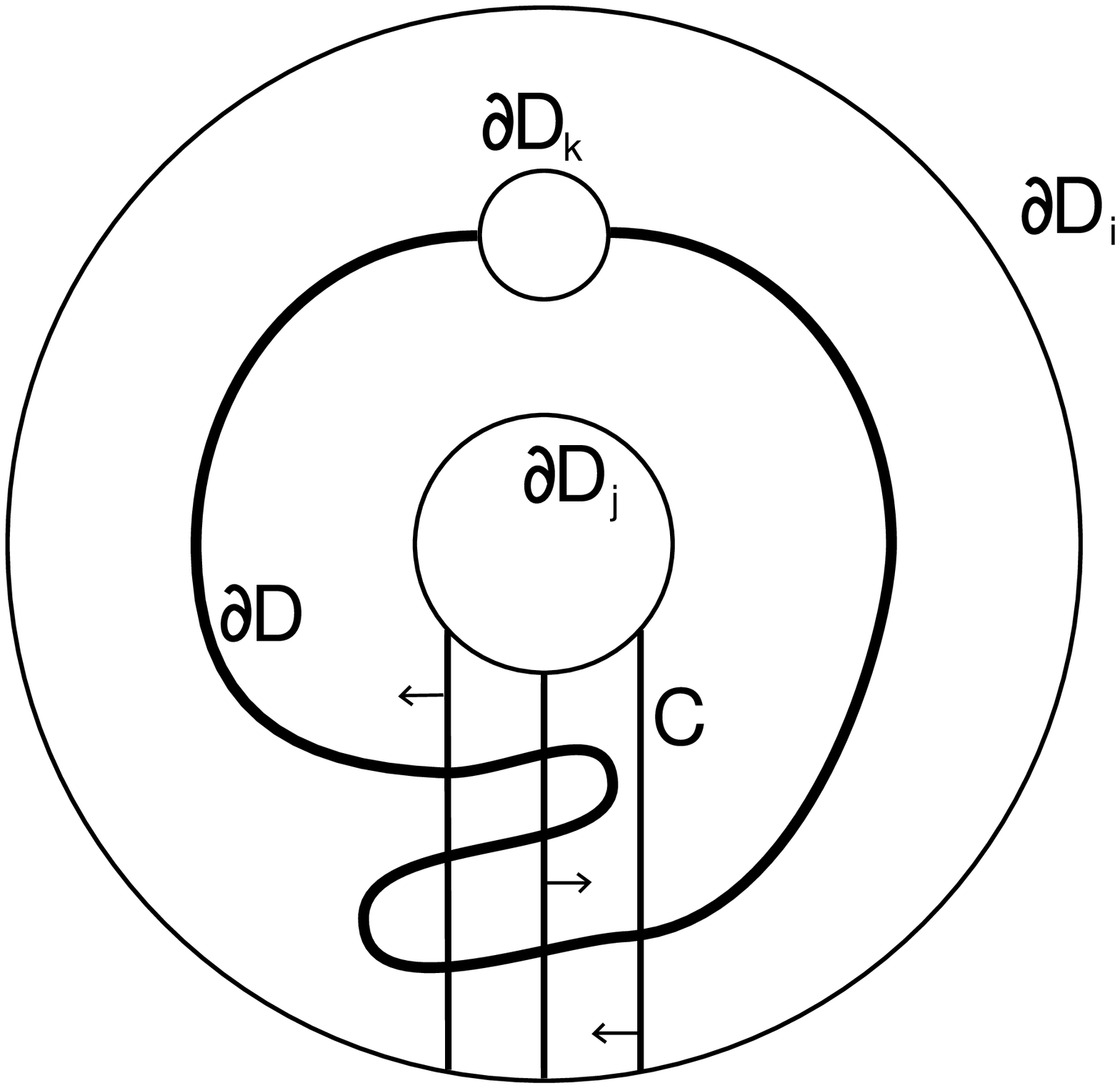}}}

\centerline{Figure 13}

\ssk

The arc $\beta$ that $\alpha$ cuts off in \del$D$, lying in one of the pairs of pants
$P$, must therefore separate the 
two \del-components \del$D_i$ and \del$D_j$ of $P$ 
which it doesn't meet. It therefore must intersect the three arcs running from 
\del$D_i$ to itself, just before and after passing through \del$D_j$, which were given 
by our hypothesis (Figure 13). As before, we may assume that $\beta$ meets 
each arc of $C$ running 
between \del$D_i$ and \del$D_j$ exactly once, since otherwise we can find a trivial
subarc of $\beta$ in $P|C$, allowing us, as before, to either reduce the number of
points of intersection of $\beta$ with $C$, or find a trivial suture after 
decomposing along $D$. But then by truncating the three arcs given by our hypothesis,
by removing the short subarcs lying at the ends between \del$D_i$ and \del$D$, we
obtain three parallel arcs whose ends all lie on the same side of \del$D$. Together 
with the obvious arcs in \del$D$, they bound a rectangle $R$ in \del$H$. 

These arcs in $C$
may not lie in \del$H|$\del$D$ (Figure 14); but since we may, as above, assume 
that every other
arc of $R\cap$\del$D$ has no trivial intersections with our triple of arcs,
some subrectangle bounded by arcs of \del$D$ will lie in \del$H|$\del$D$, with opposite
transverse orientations on the ends. The intersection of this subrectangle 
with $C$ will give us a triple of arcs with all of their ends on the
same side of \del$D$, giving us the triple of arcs in $C$ which we need to 
apply Proposition 1. Therefore decomposing $(H,C)$ along $D$ will yield a trivial 
suture, so $(H,C)$ is not disk decomposable along $D$.

\ssk

\leavevmode

\epsfxsize=2in
\centerline{{\epsfbox{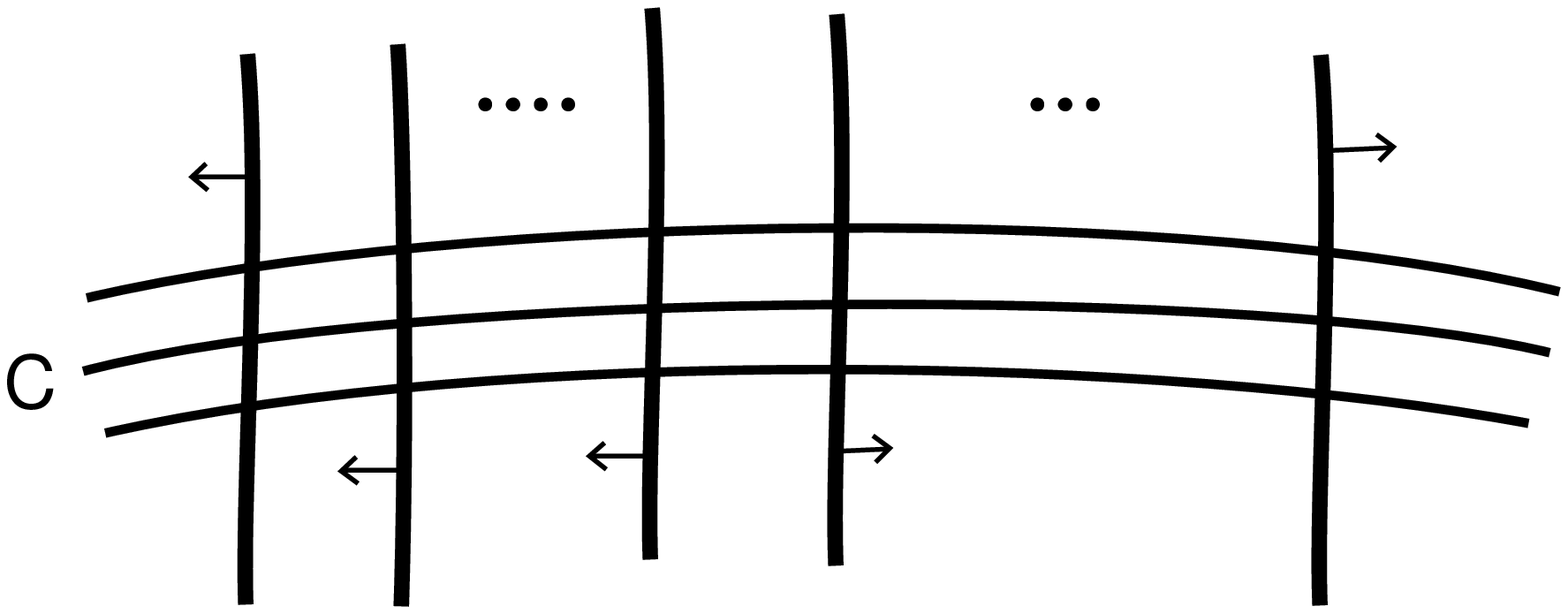}}}

\centerline{Figure 14}

\ssk

{\bf Note:} We can weaken our hypotheses somewhat while still retaining the 
conclusion. From the proof we see the we need a trio of arcs which {\it either}
end at $D_i$ after passing through $D_j$ {\it or} end at $D_j$ after passing 
through $D_i$, since we really only need the fact that the ends of the arcs are
passing {\it between} $D_i$ and $D_j$. This gives us only half as many conditions 
to check.

\heading{\S 4 \\ The examples}\endheading

\ssk

\leavevmode

\epsfxsize=3.5in
\centerline{{\epsfbox{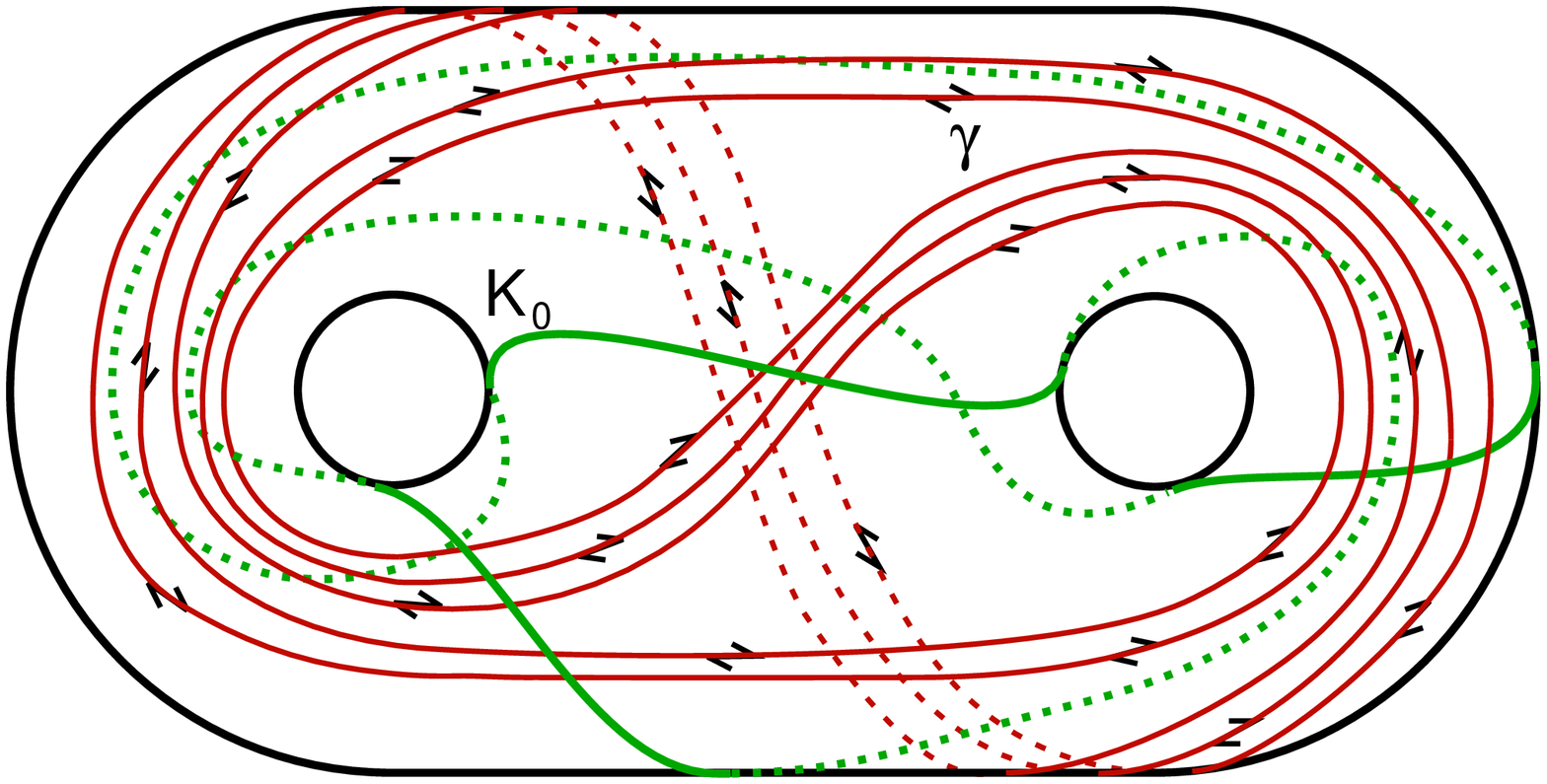}}}

\centerline{Figure 15}

\ssk

It is fairly easy to build examples of knots $K$ satisfying the conditions
of Propositions 1 and 3, by our initial Dehn twisting construction. We should 
note that the example given in Figure 7 does {\bf not} satisfy the conditions
of Proposition 3; there is no trio of arcs running from the middle disk which 
immediately run through the right hand disk on both ends. However, a still
more complicated choice of initial twisting curve $L$ will produce the examples 
we seek. Essentially, we need only make sure to choose a loop $L$ so that, for 
every choice of a pair of cutting disks, there is {\bf one} such arc
in $L$; then the fact that Dehn twisting along $L$ adds many parallel copies of
$L$ to $K_0$ will provide may parallel copies of each arc. One such example
is given in Figure 15. It is easy to verify that for each choice of
disk $D_i$, side of \del$D_i$, and choice of disk $D_j$,  $j\neq i$, there is an
arc in $L$ beginning and ending at $D_i$ on the chosen side, which immediately
passes through $D_j$ at each end (or vice versa, which suffices for our purposes
by the comment following the proof of Proposition 3). Properly chosen
subarcs of the pair of arcs shown in Figures 16ab will suffice.

\ssk

\leavevmode

\epsfxsize=3.5in
\centerline{{\epsfbox{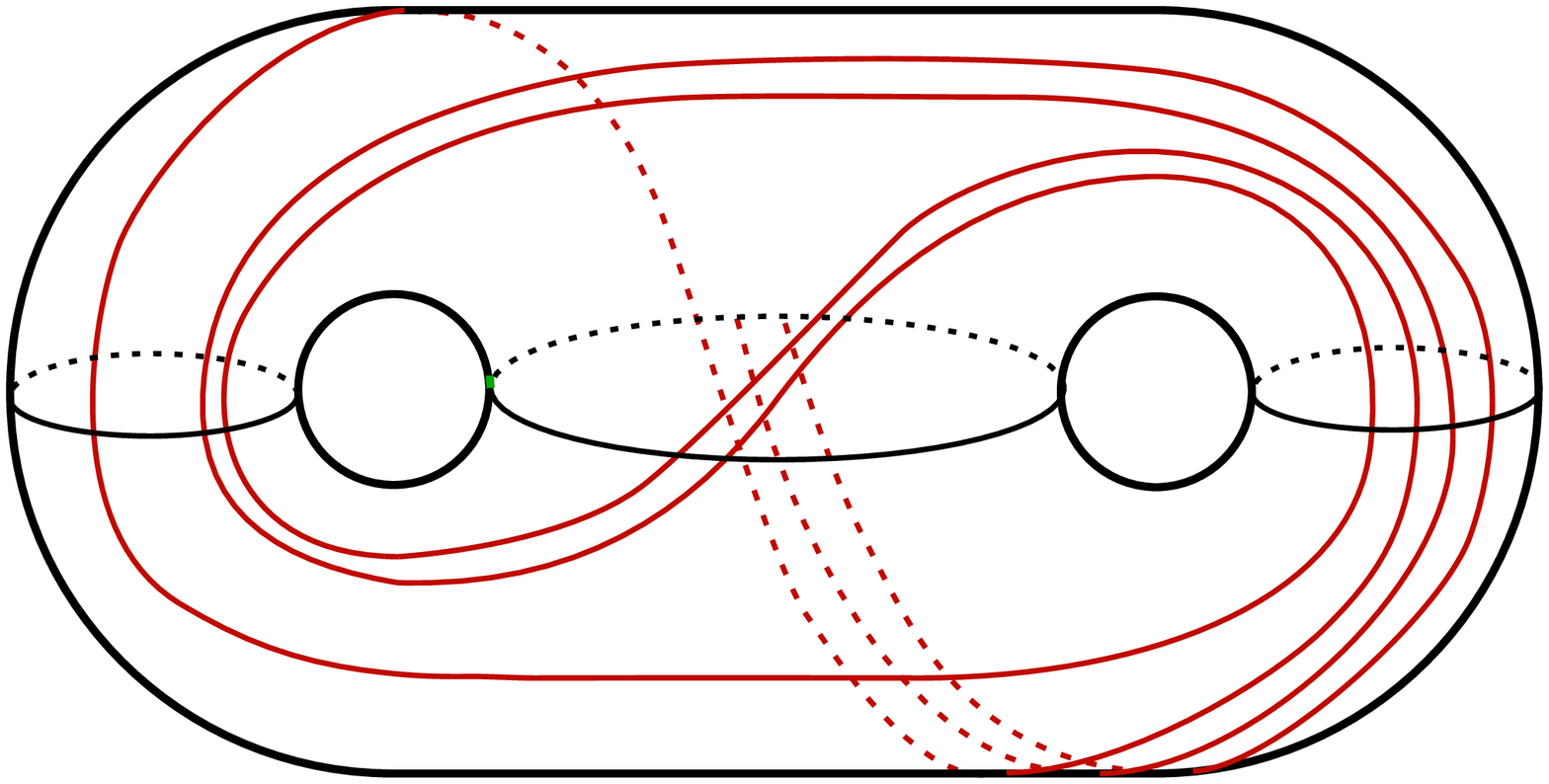}}}

\centerline{Figure 16a}

\ssk

\ssk

\leavevmode

\epsfxsize=3.5in
\centerline{{\epsfbox{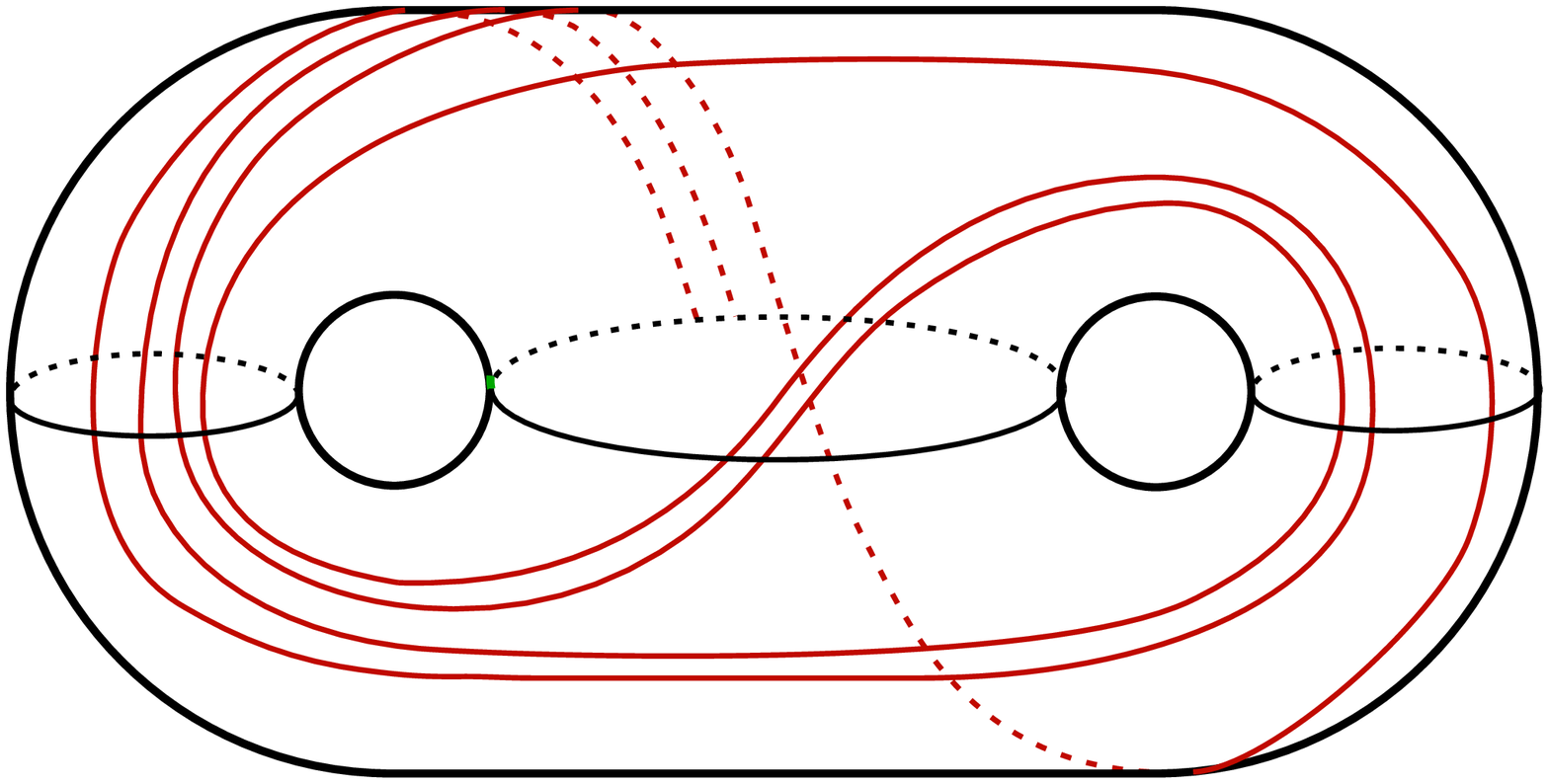}}}

\centerline{Figure 16b}

\ssk

To be certain that, when we perform a Dehn twist along $L$, the resulting loop
$C^\prime$ will have at least three arcs parallel to each of the arcs given in the
above figures, we must check that the {\it complement} of each arc $\alpha$ in $L$ 
meets 
$C$ at least three times. This is because as we traverse $\alpha$, every time we 
cross $C$ one of the arcs in $C^\prime$ parallel to $\alpha$ has been grafted 
to $C$ and (we must assume) no longer runs parallel to $\alpha$. 
Since we start with $|C\cap L|$ (= 22, in this case) arcs of $C^\prime$ running
parallel to $\alpha$ at the start, and lose one at each crossing, we simply need 
to insure that we cross $C$ no more than 19 times to ensure that three arcs
will run parallel to $\alpha$ in $C^\prime$. The reader can readily verify that
for the arcs shown in Figure 16, the complementary arcs always meet $C$ at
least 6 times, by comparing with Figure 15.

\ssk

\leavevmode

\epsfxsize=3,5in
\centerline{{\epsfbox{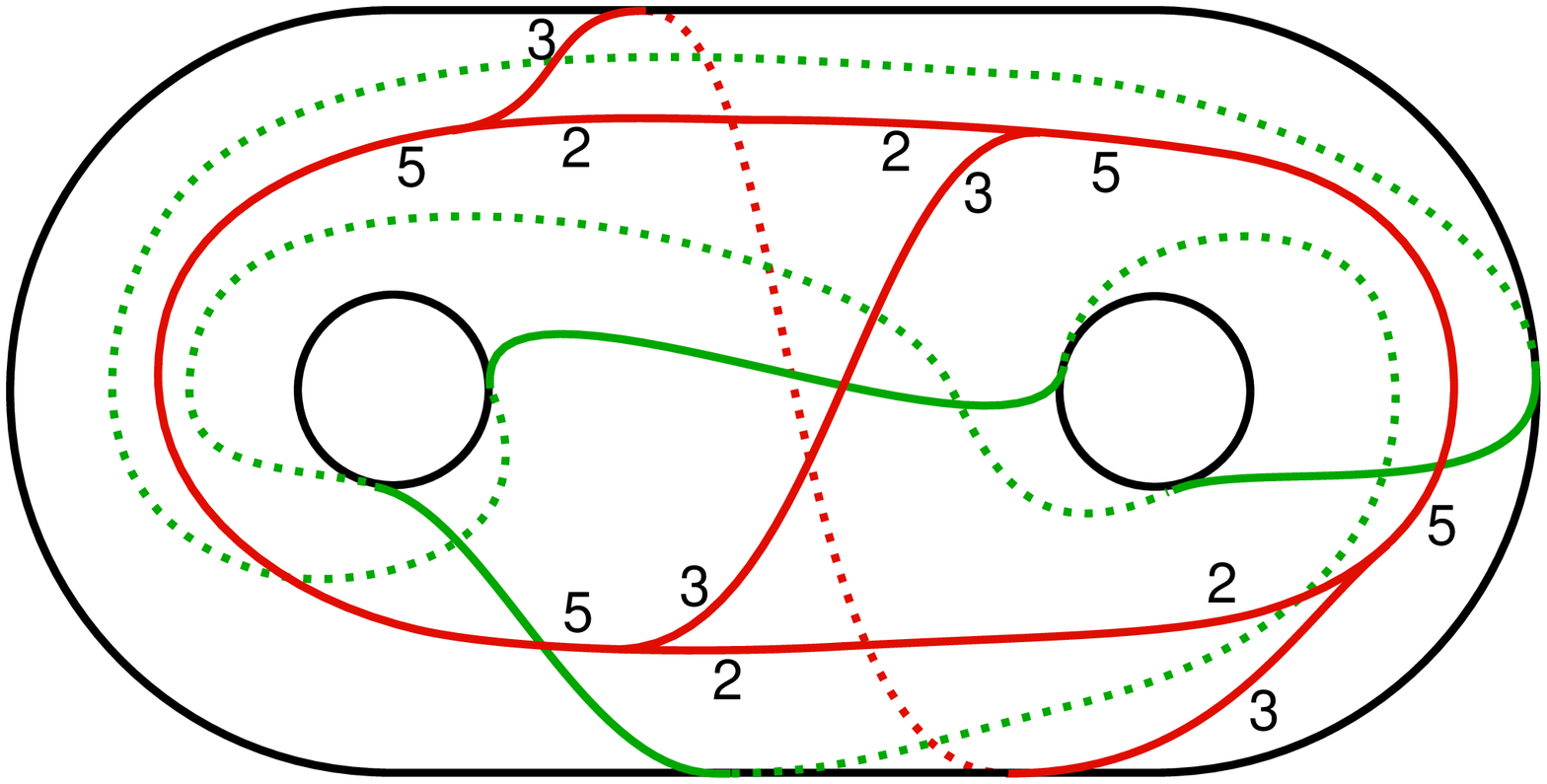}}}

\centerline{Figure 17a}

\ssk

\leavevmode

\epsfxsize=3.5in
\centerline{{\epsfbox{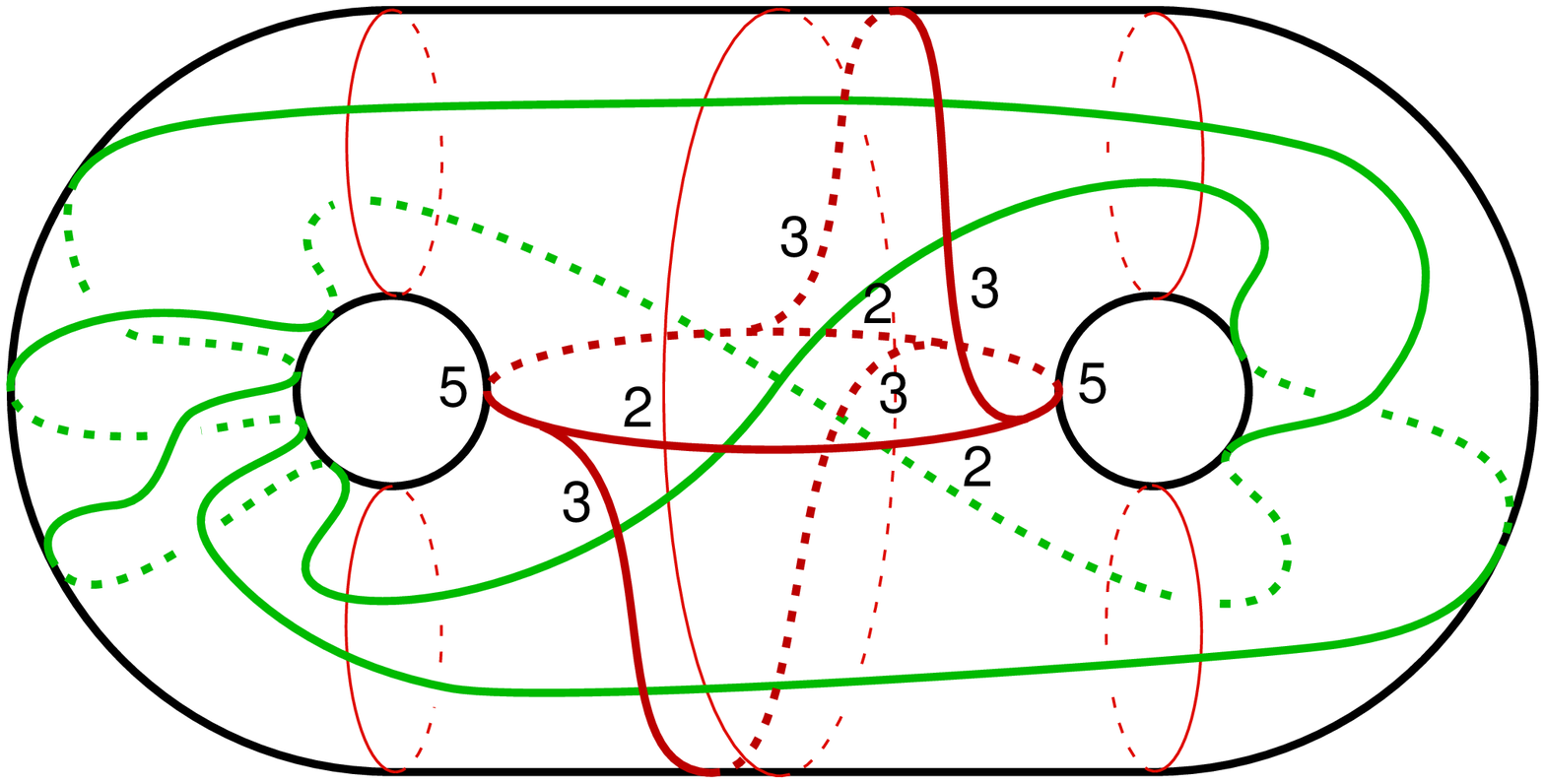}}}

\centerline{Figure 17b}

\ssk

To see what $L$ looks like in the complement of our original Seifert
surface $\Sigma$, we work with the train track $\tau$ in \del$H$ of Figure 8.
$L$ is carried by $\tau$ with weights 2, 3, and 5, as in Figure 17a. By keeping 
track of the intersections of $\tau$ with our cutting disks and $K_0$, we
can reconstruct how $\tau$ would look in the interior version of our picture
of $X_{F_0}$; see Figure 17b. This in turn allows us to reconstruct $L$, as it
sits on our 4-punctured sphere $P$ (Figure 18).

\ssk

\leavevmode

\epsfxsize=3.5in
\centerline{{\epsfbox{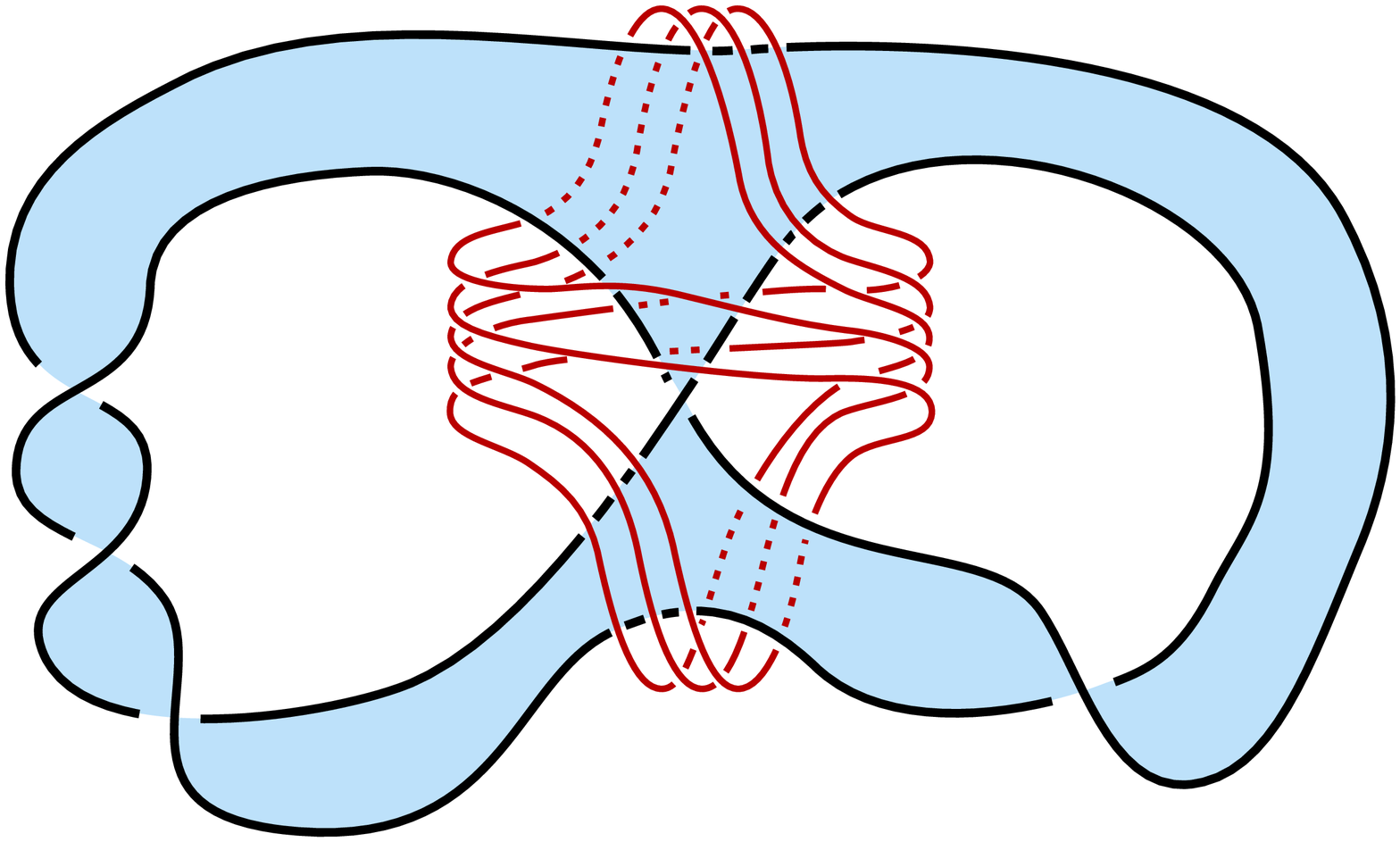}}}

\centerline{Figure 18}

\ssk

According to the computer program SnapPea [We], the knot $K$ that we obtain
from $K_0$ by 1/1 Dehn filling along $L$ is hyperbolic; by the above work,
the Seifert surface $F_0$ is carried under the Dehn filling to a 
genus one free
Seifert surface $F$ for $K$, which is not disk decomposable. 

We can 
readily construct many more such examples, since any collection of 
{\it larger} weights on the train track of Figure 16b (which represent
a connected loop, which essentially means that our replacements for 2 and 3
must be relatively prime) will also yield a knot and Seifert surface
satisfying our theorem. Similarly, $1/n$ Dehn filling along $L$ or
these more complicated loops will also suffice, since more twists 
simply provide more parallel arcs for our arguments to use. We can also add
full twists to the `arms' of our original Seifert surface, without
changing the essential features of the construction.

\heading{\S 4 \\ Concluding remarks}\endheading

The examples the we have obtained here, in some sense, manage to raise more
questions than they answer. Perhaps the most pressing question raised
is: do these knots that we build possess {\it other} Seifert surfaces
which {\it are} disk decomposable? These other surfaces must, of
course, also be free and have genus one. More generally, we might
ask:

\proclaim{Question 1} If the genus of $K$ equals the free genus of
$K$, does $K$ always possess a disk decomposable Seifert surface?
\endproclaim

One way to show that the answer to this question is `No' would be to show
that some of our examples possess only one minimal genus (free)
Seifert surface. There are several techniques for showing that a
knot possesses a unique minimal genus Seifert surface (see, e.g., 
[Ko2],[Ko3],[KK]). Most of these can be phrased as saying that the knot 
$K$ is `simple enough'; since our approach to non-disk-decomposability
is that the suture (= the knot) is 
`complicated enough', applying such techniques will no doubt
require some finesse.

While the Seifert surfaces that we build fail to be disk decomposable,
they do have minimal genus, and so Gabai [Ga1] assures us that there
is some sequence of decomposing surfaces which will split our sutured
handlebody to trivially sutured 3-balls. What we have really shown here is
that the first surface cannot be a disk. Since the decomposing surfaces
must be incompressible, they will always (inductively) split our
sutured handlebody at each stage to another sutured handlebody. The
first splitting, then, {\it cannot reduce} the genus of the 
sutured handlebody (and, except for the case of a non-separating
annulus, must {\it raise} it). An intersting question to ask, then , is:
how high must the genus of the handlebody go? Are there, for example, 
(free) Seifert surfaces (of minimal genus) for which the \underbar{first}
decomposing surface must raise the genus by an arbitrarily large
amount?

Finally, we could attempt to strengthen our result by trying to 
replace `disk decomposability' with something {\it weaker}. For
example, a disk decomposable Seifert surface is always the leaf of a 
depth one foliation of the knot exterior [Ga2], and so the knot $K$ 
must have depth [CC] (at most) one.
So one can ask the {\it weaker} question:

\proclaim{Question 2} If genus($K$) = free genus($K$), then does
$K$ have depth (at most) one?\endproclaim

An answer of `No' would be a stronger result. There is in fact a fairly
simple necessary condition for a Seifert surface to be the leaf of a
depth one foliation [CC]: the result of attaching a 2-handle to the 
suture
of the associated sutured manifold must be the 
total space of a fiber bundle over the circle. See [Ko3] for an 
example which uses the Thurston norm to check this condition.
Examples giving a negative answer to Question 1, which failed to 
satisfy this property, would also give a negative answer to Question 2.

\ssk

Nowhere in our arguments is it really essential that our system of cutting disks 
consists of {\it three} disks. The exact same conditions used here,
describing how the suture
$C$ meets a complete system of cutting disks for a higher genus handlebody, can 
therefore
be used to find higher genus examples of sutured handlebodies which are
not disk decomposable. One must use different arguments to show that the
sutured handlebody is in fact taut; the conditions we impose only guarantee that the
complement of the suture in \del$H$ is incompressible in $H$, and do not imply 
minimal genus.


\medskip
\Refs

\refstyle{A}
\widestnumber\key{MCS}

\ref\key Br
\by M. Brittenham
\paper Free genus one knots with large volume
\paperinfo preprint
\endref

\ref\key CC
\by J. Cantwell and L. Conlon
\paper Depth of knots
\jour Topology Appl
\vol 42
\yr 1991
\pages 277-289
\endref

\ref\key Ep
\by D.B.A. Epstein
\paper Curves on 2-manifolds and isotopies
\jour Acta Math
\vol 115
\yr 1966
\pages 83-107
\endref

\ref\key Ga1
\by D. Gabai
\paper Foliations and the topology of 3-manifolds
\jour J Diff Geom
\vol 18
\yr 1983
\pages
\endref 445-503

\ref\key Ga2
\bysame
\paper Foliations and genera of links
\jour Topology
\vol 23
\yr 1984
\pages 381-394
\endref

\ref\key Go
\by H. Goda
\paper A construction of taut sutured handlebodies which are not disk decomposable
\jour Kobe J Math
\vol 11
\yr 1994
\pages 107-116
\endref


\ref\key Ko1
\by T. Kobayashi
\paper Heights of simple loops and pseudo-Anosov homeomorphisms
\jour Contemp. Math.
\vol 78
\yr 1988
\pages 327-338
\endref

\ref\key Ko2
\bysame
\paper Uniqueness of minimal genus Seifert surfaces for links
\jour Topology Appl
\vol 33
\yr 1989
\pages 265-279
\endref

\ref\key Ko3
\bysame
\paper Example of hyperbolic knot which do not admit depth 1 foliation
\jour Kobe J Math 
\vol 13
\yr 1996
\pages 209-221
\endref

\ref\key KK
\by M. Kobayashi and T. Kobayashi
\paper On canonical genus and free genus of a knot
\jour J. Knot Thy. Ram.
\vol 5
\yr 1996
\pages 77-85
\endref

\ref\key Ro
\by D. Rolfsen 
\book Knots and Links 
\bookinfo Publish or Perish Press
\yr 1976
\endref

\ref\key Se
\by H. Seifert
\paper \"Uber das Geschlecht von Knoten 
\jour Math. Annalen
\vol 110
\yr 1934
\pages 571-592
\endref

\ref\key St
\by E. Starr
\paper Curves in handlebodies
\paperinfo Thesis, U Cal Berkeley, 1992
\endref

\ref\key We
\by J. Weeks
\paper SnapPea, a program for creating and studying hyperbolic 3-manifolds
\paperinfo available for download from www.geom.umn.edu
\endref

\endRefs
\enddocument